\documentclass[reqno]{amsart}
\usepackage{hyperref}
\usepackage{graphicx}
\usepackage{amssymb}

\begin{document}
    \title[\hfilneg ]
    { unique positive solution for nonlinear Caputo-type fractional $q$-difference equations with nonlocal and Stieltjes integral boundary conditions}
 \author[\hfil\hfilneg]{Ahmad Y. A. Salamooni, D. D. Pawar }
 \address{Ahmad Y. A. Salamooni \newline
     School of Mathematical Sciences, Swami Ramanand Teerth Marathwada University, Nanded-431606, India}
      \email{ayousss83@gmail.com}

 \address{D. D. Pawar \newline
    School of Mathematical Sciences, Swami Ramanand Teerth Marathwada University, Nanded-431606, India}
     \email{dypawar@yahoo.com}


    \begin{abstract}
    This paper contain a new discussion for the type of  generalized nonlinear Caputo fractional $q$-difference equations with $m$-point
    boundary value problem and Riemann-Stieltjes integral $\tilde{\alpha}[x]:=\int_{0}^{1}~x(t)d\Lambda(t).$ By applying the fixed point theorem in cones, we investigate an existence of a unique positive solution depends on $\lambda>0.$ We present some useful properties related to the Green's function
    for $m-$point boundary value problem.\\[2mm]
      \end{abstract}

\maketitle \numberwithin{equation}{section}
\newtheorem{theorem}{Theorem}[section]
\newtheorem{lemma}[theorem]{Lemma}
\newtheorem{definition}[theorem]{Definition}
\newtheorem{example}[theorem]{Example}

\newtheorem{remark}[theorem]{Remark}
\allowdisplaybreaks

 \section{Introduction}
In this paper, we have discussed the existence of unique positive solution for the $m$-point nonlinear boundary value problem
for generalized Caputo-type fractional $q$-difference equations of the
form
\begin{equation}
\left\{\begin{matrix} _{C}D^{\alpha}_{q}x(t)+\lambda h(t)f(t,x(t))=0,~~ t\in(0,1),~~n-1<\alpha\leq n,~~~ n>2,\\\\x(0)=\sum_{i=1}^{m-2}~\gamma_{i}~x(\zeta_{i}),\quad\quad\quad\\\\
 _{C}D^{2}_{q}x(0)=~_{C}D^{3}_{q}x(0)=...=~_{C}D^{n-1}_{q}x(0)=0,\quad\quad\quad\\\\
 \nu ~_{C}D_{q}x(1)-\mu\tilde{\alpha}[x]=\sum_{i=1}^{m-2} \beta_{i}~ _{C}D_{q}x(\zeta_{i}),\quad\quad\quad\end{matrix}\right.
\end{equation}
Where $_{C}D^{\alpha}_{q}$ is the generalized Caputo fractional $q$-derivative of order $\alpha,$ with \\$0<q<1,~\sum_{i=1}^{m-2}\gamma_{i}<1,~\beta_{i}\geq0,~\zeta_{i}\in(0,1),~
i=1,2,...,m-2,~\zeta_{1}<\zeta_{2}<...<\zeta_{m-2},~\lambda>0$ is a parameter,$~\nu,\mu>0,~h:[0,1]\rightarrow\mathbb{R^{+}},f:(0,1)\times\mathbb{R}\rightarrow\mathbb{R}$ is a continuous function, and $\tilde{\alpha}[x]:=\int_{0}^{1}x(t)d\Lambda(t)$ is the Riemann-Stieltjes integral with respect to the bounded variation function $\Lambda:[0,1]\rightarrow\mathbb{R}.$
\par The  fractional q-difference equations introduced and studied by Jackson [1], Adams [2], AlSalam [3] and Agarwal [4]. Due to its applicability in mathematical modeling in different branches like technical sciences, engineering, physics and biomathematics, it has  drawn wide attention to many researchers.
\par Recent studies on fractional calculus and  boundary value problems of fractional q-difference equations indicate that it is an important topic of current research. There have appeared various articles covering these problems related to finding the existence positive or nontrivial solutions for different kinds of boundary value problems such as nonlocal, integral, multiple-point, sub-strip boundary problems and some others, see [5-15] and references therein.
\par In [5], Ahmad et al. considered the following nonlinear fractional q-difference equations with nonlocal and sub-strip type boundary conditions
\[\left\{\begin{matrix} _{C}D^{\alpha}_{q}x(t)=f(t,x(t)),~~ t\in[0,1],~~1<\alpha\leq 2,~0<q<1,
\\x(0)=x_{0}+h(x),~~x(\zeta)=b\int_{\eta}^{1}x(\tau)d_{q}\tau,~~0<\zeta<\eta<1,\end{matrix}\right.\]
where $_{C}D^{\alpha}_{q}$ are the Caputo fractional $q-$derivative of order $\alpha,f:[0,1]\times\mathbb{R}\rightarrow\mathbb{R}~and~h:C([0,1],\mathbb{R})\rightarrow\mathbb{R}$ are given continuous functions, and $b$ is a real constant.
By applying Banach$'$s contraction principle and a fixed point theorem due to O$'$Regan, the existence results of the solutions were obtained.
\par The  aim of this paper is to establish the existence of unique positive solutions for the class of generalized Caputo-type fractional $q$-difference equations with $m$-point boundary value problem and Riemann-Stieltjes integral conditions by using the recent fixed point theorem in cones. Moreover, the example is given to illustrate our main results.
\section{Preliminaries}
In this section, we introduce some notations, definitions and lemmas for the theory of fractional $q$-calculus [4,8,9]. \par For $q\in(0,1),$ the class $[\hat{a}]_{q}$
is defined as
$$[\hat{a}]_{q}=\frac{1-q^{\hat{a}}}{1-q},~~\hat{a}\in\mathbb{R}.$$
The Pochhammer symbol with the $q$-analogue ($q$-shifted factorial) is defined by
$$(\hat{a};q)_{0}=1,~~(\hat{a};q)_{\ell}=\prod_{\jmath=0}^{\ell-1}(1-\hat{a}q^{\jmath}),~~\ell\in\mathbb{N}\cup{\infty}.$$
Moreover, the exponent $(\hat{a}-\hat{b})^{\ell}$ with the $q$-analogue is defined by
$$(\hat{a}-\hat{b})^{(0)}=1,~~(\hat{a}-\hat{b})^{(\ell)}=\prod_{\jmath=0}^{\ell-1}(\hat{a}-\hat{b}q^{\jmath}),~~\ell\in\mathbb{N},~~\hat{a},\hat{b}\in\mathbb{R}.$$  The $q$-gamma function $\Gamma_{q}(\hat{x})$ is defined by
$$\Gamma_{q}(\hat{x})=\frac{(1-q)^{(\hat{x}-1)}}{(1-q)^{\hat{x}-1}},~~\hat{x}\in\mathbb{R}\setminus\{0,-1,-2,...\}$$
and note that $\Gamma_{q}(\hat{x}+1)=[\hat{x}]_{q}\Gamma_{q}(\hat{x}).$
\\Now, Let $\varphi$ be the function defined on $[0,1]$ and $\alpha\geq0.$ Then, we have the following definitions:
\\\textbf{\ Definition 1[4].} The Riemann-Liouville type fractional $q$-integral of
order$~\alpha\geq0$
is $(I_{q}^{0}\varphi)(t)=\varphi(t)$ and
$$(I_{q}^{\alpha}\varphi)(t)=\frac{1}{\Gamma_{q}(\alpha)}
\int_{0}^{t}(t-q\tau)^{\alpha-1}\varphi(\tau)d_{q}\tau~,\quad\alpha>0,~t\in[0,1].$$
\\\textbf{\ Definition 2[10].} The Riemann-Liouville type fractional $q$-derivative of
order$~\alpha\geq0$
is $(D_{q}^{0}\varphi)(t)=\varphi(t)$ and
$$(D_{q}^{\alpha}\varphi)(t):=(D_{q}^{n}I_{q}^{n-\alpha}\varphi)(t),\quad\alpha>0,$$
where $n$ is a smallest integer greater than or equal to $\alpha.$
\\\textbf{\ Definition 3[10].} The Caputo type fractional $q$-derivative of
order$~\alpha>0$
is defined as
$$(_{C}D_{q}^{\alpha}\varphi)(t):=(I_{q}^{n-\alpha}D_{q}^{n}\varphi)(t),\quad\alpha>0,$$
where $n$ is a smallest integer greater than or equal to $\alpha.$
\\\textbf{\ Lemma 1[4].} Let $\varphi$ be the function defined on $[0,1]$ and $\alpha,\beta\geq0.$ Then, \\the $q$-fractional have the following property:
\\$(1)~~ (I_{q}^{\alpha}I_{q}^{\beta}\varphi)(t)=I_{q}^{\alpha+\beta}\varphi(t);~~(2)~~(D_{q}^{\alpha}I_{q}^{\alpha}\varphi)(t)=\varphi(t).$
\\\textbf{\ Lemma 2[4].}
If $\alpha\in\mathbb{R^{+}},~\gamma\in(-1,\infty),$ then
$$I_{q}^{\beta}(t)^{(\gamma)}=\frac{\Gamma_{q}(\gamma+1)}{\Gamma_{q}(\alpha+\gamma+1)}(t)^{(\alpha+\gamma)},~~~0<t<a.$$
\\\textbf{\ Lemma 3[9].}
For $\alpha>0~~and~~\alpha\in\mathbb{R^{+}}\setminus\mathbb{N},$ the following is holds:
$$(I_{q}^{\beta}~_{C}D_{q}^{\alpha}\varphi)(t)=\varphi(t)-\sum_{i=0}^{n-1}\frac{(D_{q}^{\alpha}\varphi)(0)}{\Gamma_{q}(i+1)}~~t^{i},$$
where $n$ is a smallest integer greater than or equal to $\alpha.$
\par At first, we consider the generalized Caputo type fractional $q$-difference with the following boundary value problem
\begin{equation}
\left\{\begin{matrix} _{C}D^{\alpha}_{q}x(t)+y(t)=0,~~ t\in(0,1),~~n-1<\alpha\leq n,~~~ n>2,\\\\x(0)=\sum_{i=1}^{m-2}~\gamma_{i}~x(\zeta_{i}),\quad\quad\quad\\\\
 _{C}D^{2}_{q}x(0)=~_{C}D^{3}_{q}x(0)=...=~_{C}D^{n-1}_{q}x(0)=0,\quad\quad\quad\\\\
 \nu ~_{C}D_{q}x(1)-\mu\int_{0}^{1}x(t)d\Lambda(t)=\sum_{i=1}^{m-2} \beta_{i}~ _{C}D_{q}x(\zeta_{i}),\quad\quad\quad\end{matrix}\right.
\end{equation}
and we need the following supposition:
\\$(\mathcal{H}_{1})$ The bounded variation function is $\Lambda:[0,1]\rightarrow\mathbb{R}.$ And
$$\rho=\nu-\sum_{i=1}^{m-2} \beta_{i}-\mu\mathbf{B}>0,~\mathbf{B}:=\int_{0}^{1}\big[t+\frac{\sum_{i=1}^{m-2}~\gamma_{i}\zeta_{i}}{\delta}\big]d\Lambda(t)\geq0,
~\delta=1-\sum_{i=1}^{m-2}~\gamma_{i},$$
$$\phi(\tau):=\int_{0}^{1}\big[H_{1}(t,q\tau)+H_{2}(t,q\tau;\zeta_{i})\big]d\Lambda(t)\geq0$$
for
\begin{equation}
H_{1}(t,q\tau)=\left\{\begin{matrix}\frac{\big[t+\frac{\sum_{i=1}^{m-2}~\gamma_{i}\zeta_{i}}{\delta}\big]}{2\Gamma_{q}(\alpha-1)}(1-q\tau)^{(\alpha-2)}-
\frac{(t-q\tau)^{(\alpha-1)}}{\Gamma_{q}(\alpha)},~0\leq q\tau\leq t\leq 1,\\\\
\frac{\big[t+\frac{\sum_{i=1}^{m-2}~\gamma_{i}\zeta_{i}}{\delta}\big]}{2\Gamma_{q}(\alpha-1)}(1-q\tau)^{(\alpha-2)},\quad\quad\quad\quad\quad\quad~0\leq t\leq q\tau\leq 1,
\end{matrix}\right.
\end{equation}
\begin{equation}
H_{2}(t,q\tau;\zeta_{i})=\left\{\begin{matrix}\frac{\big[t+\frac{\sum_{i=1}^{m-2}~\gamma_{i}\zeta_{i}}{\delta}\big]}{2\Gamma_{q}(\alpha-1)}(1-q\tau)^{(\alpha-2)}-
\frac{\sum_{i=1}^{m-2}~\gamma_{i}}{\delta\Gamma_{q}(\alpha)}(\zeta_{i}-q\tau)^{(\alpha-1)},~0\leq q\tau\leq\zeta_{i}\leq 1,\\\\
\frac{\big[t+\frac{\sum_{i=1}^{m-2}~\gamma_{i}\zeta_{i}}{\delta}\big]}{2\Gamma_{q}(\alpha-1)}(1-q\tau)^{(\alpha-2)},
\quad\quad\quad\quad\quad\quad\quad\quad\quad\quad\quad\quad~0\leq\zeta_{i}\leq q\tau\leq 1.
\end{matrix}\right.
\end{equation}
And let
\begin{equation}
H_{3}(t,q\tau;\zeta_{i})=\frac{\big[t+\frac{\sum_{i=1}^{m-2}~\gamma_{i}\zeta_{i}}{\delta}\big]}{\Gamma_{q}(\alpha-1)}\left\{\begin{matrix}(1-q\tau)^{(\alpha-2)}-
(\zeta_{i}-q\tau)^{(\alpha-2)},~0\leq q\tau\leq\zeta_{i}\leq 1,\\\\(1-q\tau)^{(\alpha-2)},
\quad\quad\quad\quad\quad\quad\quad\quad~0\leq\zeta_{i}\leq q\tau\leq 1.
\end{matrix}\right.
\end{equation}
\\\textbf{\ Lemma 4.}
Let $y\in\mathcal{C}[0,1]$ and assume that $(\mathcal{H}_{1})$ holds, then the boundary value problem $(2.1)$ has the unique solution $x$ given in the form
\begin{equation}
x(t)=\int_{0}^{1}~G(t,q\tau;\zeta_{i})~y(\tau)~d_{q}\tau,
\end{equation}
where
\begin{equation}
 \begin{gathered}
G(t,q\tau;\zeta_{i})=H_{1}(t,q\tau)+H_{2}(t,q\tau;\zeta_{i})\quad\quad\quad\quad\quad\quad\quad\quad\quad\quad\quad\\+\frac{\sum_{i=1}^{m-2}~ \beta_{i}}{\rho}~~H_{3}(t,q\tau;\zeta_{i})\quad\quad\quad\\+\big[t+\frac{\sum_{i=1}^{m-2}~\gamma_{i}\zeta_{i}}{\delta}\big]~~\frac{\mu}{\rho}~~\phi(\tau).
\end{gathered}
\end{equation}
\\ \textbf{\ Proof.} In view of the Lemma 3, the solution of
the generalized Caputo type fractional $q$-difference equation $(2.1)$ can be written as
\begin{equation}
 \begin{gathered}
x(t)=-\frac{1}{\Gamma_{q}(\alpha)}\int_{0}^{1}~(t-q\tau)^{(\alpha-1)}~y(\tau)~d_{q}\tau
+c_{0}+c_{1}t+...+c_{n-1}t^{n-1},\\
_{C}D_{q}x(t)=-\frac{1}{\Gamma_{q}(\alpha-1)}\int_{0}^{1}~(t-q\tau)^{(\alpha-2)}~y(\tau)~d_{q}\tau
+c_{1}\quad\quad\quad\quad\quad\quad\quad\quad\\ \quad\quad+2c_{2}t+3c_{3}t^{2}+...+(n-1)c_{n-1}t^{n-2},\\
\vdots\\
_{C}D_{q}^{(n-1)}x(t)=-\frac{1}{\Gamma_{q}(\alpha-n+1)}\int_{0}^{1}~(t-q\tau)^{(\alpha-n)}~y(\tau)~d_{q}\tau
+(n-1)!~c_{n-1},\\
\end{gathered}
\end{equation}
where $c_{0}, c_{1}, ..., c_{n-1}\in\mathbb{R}$ are arbitrary constants. Applying
the boundary conditions, we found that
\begin{equation}\left\{\begin{matrix}
 \begin{gathered}
 c_{2}= c_{3}= c_{4}= ...= c_{n-1}=0,\quad\quad\quad\quad\quad\quad\quad\quad\\
c_{1}=I_{q}^{\alpha-1}y(1)+\frac{\sum_{i=1}^{m-2}~\beta_{i}}{\rho}~\big[I_{q}^{\alpha-1}y(1)-I_{q}^{\alpha-1}y(\zeta_{i})\big]
\quad\quad\quad\quad\\
\quad+\frac{\mu}{\rho}~\int_{0}^{1}\Bigg[\big[t+\frac{\sum_{i=1}^{m-2}~\gamma_{i}\zeta_{i}}{\delta}\big]~I_{q}^{\alpha-1}y(1)
\\
\quad\quad\quad\quad\quad\quad\quad\quad\quad\quad-\big[I_{q}^{\alpha}y(t)+
\frac{\sum_{i=1}^{m-2}~\gamma_{i}}{\delta}~I_{q}^{\alpha}y(\zeta_{i})\big]\Bigg]d\Lambda(t)\\
c_{0}=\frac{\sum_{i=1}^{m-2}~\gamma_{i}\zeta_{i}}{\delta}\Big[I_{q}^{\alpha-1}y(1)+~
\frac{\sum_{i=1}^{m-2}~\beta_{i}}{\rho}~\big[I_{q}^{\alpha-1}y(1)-I_{q}^{\alpha-1}y(\zeta_{i})\big]\Big]
\\+\frac{\mu}{\rho}\big[\frac{\sum_{i=1}^{m-2}~\gamma_{i}\zeta_{i}}{\delta}\big]
\int_{0}^{1}\Bigg[\big[t+\frac{\sum_{i=1}^{m-2}~\gamma_{i}\zeta_{i}}{\delta}\big]~I_{q}^{\alpha-1}y(1)\\
\quad\quad\quad\quad\quad\quad\quad\quad\quad\quad\quad\quad\quad\quad-\big[I_{q}^{\alpha}y(t)+
\frac{\sum_{i=1}^{m-2}~\gamma_{i}}{\delta}~I_{q}^{\alpha}y(\zeta_{i})\big]\Bigg]d\Lambda(t)\\
-\frac{\sum_{i=1}^{m-2}~\gamma_{i}}{\delta}I_{q}^{\alpha}y(\zeta_{i}).
\end{gathered}\end{matrix}\right.
\end{equation}
Substituting the values of $c_{0}, c_{1}, ..., c_{n-1},$ we deduce that
\begin{equation}\left\{\begin{matrix}
 \begin{gathered}
x(t)=\frac{\big[t+\frac{\sum_{i=1}^{m-2}~\gamma_{i}\zeta_{i}}{\delta}\big]}{\Gamma_{q}(\alpha-1)}
\int_{0}^{1}~(1-q\tau)^{(\alpha-2)}~y(\tau)~d_{q}\tau\quad\quad\quad\quad\quad\quad\quad\quad\quad\quad\quad\quad\quad\quad\quad\\-\frac{1}{\Gamma_{q}(\alpha)}
\int_{0}^{1}~(t-q\tau)^{(\alpha-1)}~y(\tau)~d_{q}\tau\quad\quad\quad\quad\quad\quad
\\-\frac{\sum_{i=1}^{m-2}~\gamma_{i}}{\delta\Gamma_{q}(\alpha)}\int_{0}^{\zeta_{i}}~(\zeta_{i}-q\tau)^{(\alpha-1)}~y(\tau)~d_{q}\tau
\\+\frac{\mu}{\rho}\big[t+\frac{\sum_{i=1}^{m-2}~\gamma_{i}\zeta_{i}}{\delta}\big]\Bigg[\int_{0}^{1}\int_{0}^{1}
\frac{\big[t+\frac{\sum_{i=1}^{m-2}~\gamma_{i}\zeta_{i}}{\delta}\big]}{\Gamma_{q}(\alpha-1)}(1-q\tau)^{(\alpha-2)}~y(\tau)~d_{q}\tau~d\Lambda(t)
\\ \quad\quad\quad\quad-\int_{0}^{1}\int_{0}^{t}\frac{(t-q\tau)^{(\alpha-1)}}{\Gamma_{q}(\alpha)}~y(\tau)~d_{q}\tau~d\Lambda(t)
\\ \quad\quad\quad\quad\quad\quad\quad\quad\quad\quad-\frac{\sum_{i=1}^{m-2}~\gamma_{i}}{\delta}\int_{0}^{1}\int_{0}^{\zeta_{i}}
\frac{(\zeta_{i}-q\tau)^{(\alpha-1)}}{\Gamma_{q}(\alpha)}~y(\tau)~d_{q}\tau~d\Lambda(t)\Bigg],
\end{gathered}\end{matrix}\right.
\end{equation}
that is
\begin{equation}
 \begin{gathered}
x(t)=\int_{0}^{1}\big[H_{1}(t,q\tau)+H_{2}(t,q\tau;\zeta_{i})\big]~y(\tau)~d_{q}\tau\quad\quad\quad\quad\quad\quad\quad\quad\quad\quad\quad\\
+\int_{0}^{1}\frac{\sum_{i=1}^{m-2}~ \beta_{i}}{\rho}~H_{3}(t,q\tau;\zeta_{i})~y(\tau)~d_{q}\tau\quad\quad\quad\\
+\int_{0}^{1}\big[t+\frac{\sum_{i=1}^{m-2}~\gamma_{i}\zeta_{i}}{\delta}\big]~\frac{\mu}{\rho}~\phi(\tau)~y(\tau)~d_{q}\tau.
\end{gathered}
\end{equation}
Hence, the proof of this Lemma is complete.$\quad\quad\quad\quad\quad\quad\quad\quad\quad\quad\quad\quad\quad\quad\quad\quad\Box$
\\\textbf{\ Lemma 5.} $H_{1}(t,q\tau), H_{2}(t,q\tau;\zeta_{i}), H_{3}(t,q\tau;\zeta_{i}) ~and~ G(t,q\tau;\zeta_{i})$
are the Green's functions defined by $(2.2), (2.3), (2.4) ~and~ (2.6)$ respectively, satisfy the following properties:
\begin{equation}
 \begin{gathered}
(i)~~H_{1}(t,q\tau), H_{2}(t,q\tau;\zeta_{i}), H_{3}(t,q\tau;\zeta_{i})>0,~\forall t,\tau\in(0,1),\quad\quad\quad\quad\quad\quad\quad\quad\quad\quad\quad\\\&~~
H_{1}(t,q\tau)+ H_{2}(t,q\tau;\zeta_{i})\leq~\big[t+\frac{\sum_{i=1}^{m-2}~\gamma_{i}\zeta_{i}}{\delta}\big]\frac{(1-q\tau)^{(\alpha-2)}}{\Gamma_{q}(\alpha-1)}
~\leq\frac{\sigma~(1-q\tau)^{(\alpha-2)}}{\Gamma_{q}(\alpha-1)},\\\&~~
\frac{1}{\sigma}H_{3}(t,q\tau;\zeta_{i})\leq~ H_{3}(t,q\tau;\zeta_{i})\leq~\big[t+\frac{\sum_{i=1}^{m-2}~\gamma_{i}\zeta_{i}}{\delta}\big]\frac{(1-q\tau)^{(\alpha-2)}}{\Gamma_{q}(\alpha-1)},\\ \forall t,\tau\in[0,1],~\sigma= 1+\frac{\sum_{i=1}^{m-2}~\gamma_{i}\zeta_{i}}{\delta},~i=1,2,...,m-2.
\end{gathered}
\end{equation}
\begin{equation}
 \begin{gathered}
(ii)~~\forall t,\tau\in[0,1],~i=1,2,...,m-2,~ we~have\quad\quad\quad\quad\quad\quad\quad\quad\quad\quad\quad\quad\quad\quad\quad\quad\\ \frac{\big[t+\frac{\sum_{i=1}^{m-2}~\gamma_{i}\zeta_{i}}{\delta}\big]}{\sigma}\Big[H_{1}(1,q\tau)\quad\quad\quad\quad\quad\quad\quad\quad\quad
\quad\quad\quad\quad\quad\quad\quad\quad\quad\quad\quad\quad\quad\quad\quad\\+H_{2}(1,q\tau;1)\Big]\leq~
H_{1}(t,q\tau)+H_{2}(t,q\tau;\zeta_{i}) \leq~H_{1}(1,q\tau)+H_{2}(1,q\tau;1).
\end{gathered}
\end{equation}
\begin{equation}
 \begin{gathered}
(iii)~~\frac{\big[t+\frac{\sum_{i=1}^{m-2}~\gamma_{i}\zeta_{i}}{\delta}\big]}{\sigma}\psi_{1}(\tau)\leq~
G(t,q\tau;\zeta_{i})\leq~\big[t+\frac{\sum_{i=1}^{m-2}~\gamma_{i}\zeta_{i}}{\delta}\big]\psi_{2}(\tau),\\
here,~~\psi_{1}(\tau)=\big[\frac{\sigma}{\Gamma_{q}(\alpha-1)}(1-q\tau)^{(\alpha-2)}-
\frac{1}{\delta~\Gamma_{q}(\alpha)}(1-q\tau)^{(\alpha-1)}\big]\hat{P}\\ \quad\quad\quad\quad+\big[(1-q\tau)^{(\alpha-2)}-
(\zeta_{m-2}-q\tau)^{(\alpha-2)}\big]\tilde{P},\\ \&~~\psi_{2}(\tau)=\big[\frac{1}{\Gamma_{q}(\alpha-1)}\hat{P}+\tilde{P}\big](1-q\tau)^{(\alpha-2)},
\quad\quad\quad\quad\quad\quad\quad\\~~where\quad
\hat{P}=1+\frac{\mu~\mathbf{B}}{\rho},~~
\tilde{P}=\frac{\sum_{i=1}^{m-2} \beta_{i}}{\rho~\Gamma_{q}(\alpha-1)},~~\forall t,\tau\in[0,1].
\end{gathered}
\end{equation}
\\ \textbf{\ Proof.}  Easy to see that $(i)\&(ii)$ are holds.
\\Now, we prove $(iii)$ with the help of $(i)\&(ii),$ so we have
\begin{align*}
G(t,q\tau;\zeta_{i})&=H_{1}(t,q\tau)+H_{2}(t,q\tau;\zeta_{i})+\frac{\sum_{i=1}^{m-2}~ \beta_{i}}{\rho}~~H_{3}(t,q\tau;\zeta_{i})+
\big[t+\frac{\sum_{i=1}^{m-2}~\gamma_{i}\zeta_{i}}{\delta}\big]~~\frac{\mu}{\rho}~~\phi(\tau)\\&\leq
\bigg[\frac{\big[t+\frac{\sum_{i=1}^{m-2}~\gamma_{i}\zeta_{i}}{\delta}\big]}{\Gamma_{q}(\alpha-1)}(1-q\tau)^{(\alpha-2)}\bigg]\bigg[1+\frac{\sum_{i=1}^{m-2} \beta_{i}}{\rho}\bigg]\\&\quad\quad\quad+\frac{\mu~\big[t+\frac{\sum_{i=1}^{m-2}~\gamma_{i}\zeta_{i}}{\delta}\big]}{\rho~\Gamma_{q}(\alpha-1)}
\int_{0}^{1}\big[t+\frac{\sum_{i=1}^{m-2}~\gamma_{i}\zeta_{i}}{\delta}\big](1-q\tau)^{(\alpha-2)}d\Lambda(t)\\&=
\big[t+\frac{\sum_{i=1}^{m-2}~\gamma_{i}\zeta_{i}}{\delta}\big]\big[\frac{1}{\Gamma_{q}(\alpha-1)}\hat{P}+\tilde{P}\big](1-q\tau)^{(\alpha-2)}
\end{align*}
\begin{equation}
\quad\quad=\big[t+\frac{\sum_{i=1}^{m-2}~\gamma_{i}\zeta_{i}}{\delta}\big]\psi_{2}(\tau)
\quad\quad\quad\quad\quad\quad\quad\quad\quad\quad\quad\quad\quad\quad\quad
\end{equation}
Furthermore, we have
\begin{align*}
G(t,q\tau;\zeta_{i})&\geq\frac{\big[t+\frac{\sum_{i=1}^{m-2}~\gamma_{i}\zeta_{i}}{\delta}\big]}{\sigma}
\Bigg[\frac{\sigma~(1-q\tau)^{(\alpha-2)}}{\Gamma_{q}(\alpha-1)}-\frac{(1-q\tau)^{(\alpha-1)}}{\delta~\Gamma_{q}(\alpha)}
\quad\quad\quad\quad\quad\quad\quad\quad\\&\quad\quad\quad\quad
+\frac{\sum_{i=1}^{m-2} \beta_{i}}{\rho~\Gamma_{q}(\alpha-1)}\big[(1-q\tau)^{(\alpha-2)}-(\zeta_{m-2}-q\tau)^{(\alpha-2)}\big]
\\&\quad\quad\quad\quad+\frac{\mu}{\rho}\int_{0}^{1}\big[t+\frac{\sum_{i=1}^{m-2}~\gamma_{i}\zeta_{i}}{\delta}\big]
\bigg[\frac{\sigma~(1-q\tau)^{(\alpha-2)}}{\Gamma_{q}(\alpha-1)}-\frac{(1-q\tau)^{(\alpha-1)}}{\delta~\Gamma_{q}(\alpha)}\bigg]d\Lambda(t)\Bigg]\\&=
\frac{\big[t+\frac{\sum_{i=1}^{m-2}~\gamma_{i}\zeta_{i}}{\delta}\big]}{\sigma}
\Bigg[\bigg[\frac{\sigma~(1-q\tau)^{(\alpha-2)}}{\Gamma_{q}(\alpha-1)}-\frac{(1-q\tau)^{(\alpha-1)}}{\delta~\Gamma_{q}(\alpha)}\bigg]
\bigg[1+\frac{\mu~\mathbf{B}}{\rho}\bigg]\\&\quad\quad\quad\quad+\frac{\sum_{i=1}^{m-2} \beta_{i}}{\rho~\Gamma_{q}(\alpha-1)}
\big[(1-q\tau)^{(\alpha-2)}-
(\zeta_{m-2}-q\tau)^{(\alpha-2)}\big]\Bigg]
\end{align*}
\begin{equation}
\quad\quad=\frac{\big[t+\frac{\sum_{i=1}^{m-2}~\gamma_{i}\zeta_{i}}{\delta}\big]}{\sigma}\psi_{1}(\tau)
\quad\quad\quad\quad\quad\quad\quad\quad\quad\quad\quad\quad\quad\quad\quad
\end{equation}
Hence, the proof of this Lemma is complete.$\quad\quad\quad\quad\quad\quad\quad\quad\quad\quad\quad\quad\quad\quad\quad\quad\Box$
\par Let the real Banach space $(K,\|.\|)$ is partially ordered by the cone $E\subset K,$ and we denote
a zero element of $K$ by $\phi.$ For the fixed $p>\phi(i.e.,p\geq\phi~and~p\neq\phi),$ the set $E_{p}$ defined by $E_{p}=\{x\in K| x\sim p\},$
here $\forall x,y\in K,$ the notation $ x\sim y$ means that $\exists \hat{\mu},\hat{\nu}>0$ such that $\hat{\mu}x\leq y\leq\hat{\nu}x,$ and
easy to see that $E_{p}\subset E.$
\par Now, we will give the recent fixed point theorems, which play a key role in the next analysis.
\\\textbf{\ Lemma 6[16]} Let $E$ be a normal cone in the real Banach space $K$ with $p>\phi,~and~\mathcal{T}:E\rightarrow E$ is
an increasing operator, satisfying:
\\$(i)$ there is $p_{0}\in E_{p}$ such that $\mathcal{T}p_{0}\in E_{p};$
\\$(ii)$ for any $x\in E,\ell\in(0,1),$ there exists $y(\ell)\in(\ell,1)$ such that $\mathcal{T}(\ell x)(t)\geq y(\ell)\mathcal{T}x(t).$\\Then:
\\(1) the operator equation $\mathcal{T}x=x$ has a unique solution $x^{*}$ in $E_{p};$
\\(2) for any initial value $x_{0}\in E_{p},$ constructing successively the sequence $x_{n}=\mathcal{T}x_{n-1},~n=1,2,...,$ we have $x_{n}\rightarrow x^{*}$ as $n\rightarrow\infty.$
\\\textbf{\ Remark 1.[16]}  It is remarkable that \\(1) if $y(\ell)=\ell,\ell\in(0,1)$ in Lemma 6, then the operator $\mathcal{T}:E\rightarrow E$ is said to be sub-homogeneous;\\(2) if $y(\ell)=\ell^{r},0\leq r<1$ in Lemma 6, then the operator $\mathcal{T}:E\rightarrow E$ is said to be $r-concave.$
\\\textbf{\ Lemma 7[16]} On the basis of Lemma 6, if $x_{\lambda}$ is a unique solution of the operator
equation $\mathcal{T}x=\lambda x$ for $\lambda>0,$ Then we have the following conclusions:
\\$(i)$ $x_{\lambda}$ is strictly decreasing in $\lambda,$ namely, $0<\lambda_{1}<\lambda_{2}$ implies $x\lambda_{1}>x\lambda_{2};$
\\$(ii)$ if there exists $r\in(0,1)$ such that $y(\ell)\geq\ell^{r}$ for $\ell\in(0,1),$ then $x_{\lambda}$ is continuous
in $\lambda,$ namely, $\lambda\rightarrow\lambda_{0}(\lambda_{0}>0)$ implies $\|x_{\lambda}-x_{\lambda_{0}}\|\rightarrow0;$
\\$(iii)$ $\lim_{\lambda\rightarrow+\infty} \|x_{\lambda}\|=0,~\lim_{\lambda\rightarrow0} \|x_{\lambda}\|=+\infty.$

\section{Main Results}
In this section, we present the existence of a unique positive solution for boundary value problem (1.1) by using the recent fixed point theorem in cones. Let $K=C[1,0]$ be a Banach space endowed with the norm $\|x\|=max\{|x(t)|:t\in[0,1]\}$ and let $E$ be the standard cone $E=\{x\in K:x(t)\geq0,t\in[0,1]\}.$ So $E$ is normal. The operator $\mathcal{T}:K\rightarrow K$ is given by
\begin{equation}
\mathcal{T}x(t)=\int_{0}^{1}G(t,q\tau;\zeta_{i})~h(\tau)~f(\tau,x(\tau))~d_{q}\tau,~~t\in[0,1].
\end{equation}
\\\textbf{\ Theorem 1.} Assume that $(\mathcal{H}_{1})$ holds. In addition,
\\$(\mathcal{H}_{2})~f\in C([1,0]\times\mathbb{R}^{+},\mathbb{R}^{+}),h\in C([1,0]\times\mathbb{R}^{+}),f(t,0)\geq0,~f(t,0)\not\equiv0,h(t)\not\equiv0,for~~ t\in[0,1];$
\\$(\mathcal{H}_{3})~$for each $t\in[0,1],$ the function $f(t,x)$ is increasing in $x;$
\\$(\mathcal{H}_{4})~$for any $\ell\in(0,1),$ there exists $y(\ell)\in(\ell,1)$ such that $f(t,\ell x)\geq y(\ell)f(t,x),\forall t\in[0,1],\forall x\in\mathbb{R}^{+}.$\\Then:
\\$1.$ For any fixed $\lambda>0,$ the problem (1.1) has a unique positive solution $x^{*}_{\lambda}\in E_{p},~where~p(t)=\frac{\big[t+\frac{\sum_{i=1}^{m-2}\gamma_{i}\zeta_{i}}{\delta}\big]}{\sigma},~t\in[0,1].$ In addition, for any initial value $x_{0}\in E_{p},$ constructing the sequence
\begin{equation}
x_{n}(t)=\lambda\int_{0}^{1}G(t,q\tau;\zeta_{i})~h(\tau)~f(\tau,x_{n-1}(\tau))~d_{q}\tau,~n=1,2,...,
\end{equation}
one has $\lim_{n\rightarrow+\infty}x_{n}(t)=x^{*}_{\lambda}(t),~t\in[0,1];$
\\$2.$ $x^{*}_{\lambda}$ is strictly increasing in $\lambda,$ namely, $0<\lambda_{1}<\lambda_{2}$ implies $x^{*}_{\lambda_{1}}<x^{*}_{\lambda_{2}};$
\\$3.$ If there exists $r\in(0,1)$ such that $y(\ell)\geq\ell^{r}$ for $\ell\in(0,1),$ then $x^{*}_{\lambda}$ is continuous
in $\lambda,$ that is, $\lambda\rightarrow\lambda_{0}(\lambda_{0}>0)$ means $\|x_{\lambda}-x_{\lambda_{0}}\|\rightarrow0;$
\\$4.~\lim_{\lambda\rightarrow+\infty} \|x_{\lambda}\|=+\infty,~\lim_{\lambda\rightarrow0^{+}} \|x_{\lambda}\|=0.$
\\ \textbf{\ Proof.} According to the Lemma 4, it is easy to see that $x(t)$ is a solution of problem (1.1) if, and only if, $x(t)=\lambda \mathcal{T}x(t).$
By $(\mathcal{H}_{2})$ and Lemma 5, $T:E\rightarrow E$ is clear. From the assumption $(\mathcal{H}_{3})~$ we can easily obtain $\mathcal{T}:E\rightarrow E$ is increasing.
\par Now, we will check that $\mathcal{T}$ satisfies all the assumptions of Lemma 6. For a condition $(i)$ of Lemma 6, we put $p(t)=\frac{\big[t+\frac{\sum_{i=1}^{m-2}\gamma_{i}\zeta_{i}}{\delta}\big]}{\sigma}\geq0,$ that is, $p\in E,$ and we will show that $\mathcal{T}p\in E_{p}.$
 Set
\begin{equation}
\theta_{1}=\int_{0}^{1}\psi_{1}(\tau)~h(\tau)~f(\tau,0)~d_{q}\tau,\quad\theta_{2}=\int_{0}^{1}\sigma~\psi_{2}(\tau)~h(\tau)~f(\tau,1)~d_{q}\tau,
\end{equation}
since $\nu\geq0,~f$ is increasing with $f(t,0)\geq0,~f(t,0)\not\equiv0,h(t)\not\equiv0,$ in view of Lemma 5, we can easily obtain $0<\theta_{1}\leq\theta_{2}.$
 By $(\mathcal{H}_{3}),$ we have
\begin{align*}
\mathcal{T}p(t)&=\int_{0}^{1}G(t,q\tau;\zeta_{i})~h(\tau)~f(\tau,\frac{\big[t+\frac{\sum_{i=1}^{m-2}\gamma_{i}\zeta_{i}}{\delta}\big]}{\sigma})~d_{q}\tau
\quad\quad\quad\quad\quad\quad\\&
\leq\int_{0}^{1}\big[t+\frac{\sum_{i=1}^{m-2}\gamma_{i}\zeta_{i}}{\delta}\big]\psi_{2}(\tau)h(\tau)~f(\tau,1)~d_{q}\tau\quad\quad\quad\quad\quad\quad
\end{align*}
\begin{equation}
~=\int_{0}^{1}\sigma~\psi_{2}(\tau)~h(\tau)~f(\tau,1)~d_{q}\tau\cdot\frac{\big[t+\frac{\sum_{i=1}^{m-2}\gamma_{i}\zeta_{i}}{\delta}\big]}{\sigma}=\theta_{2}p(t).
\end{equation}
Furthermore, we get
\begin{align*}
\mathcal{T}p(t)\geq\int_{0}^{1}\frac{\big[t+\frac{\sum_{i=1}^{m-2}\gamma_{i}\zeta_{i}}{\delta}\big]}{\sigma}\psi_{1}(\tau)h(\tau)~f(\tau,0)~d_{q}\tau
\quad\quad\quad\quad\quad\quad\quad\quad
\end{align*}
\begin{equation}
=\int_{0}^{1}\psi_{1}(\tau)h(\tau)~f(\tau,0)~d_{q}\tau\cdot\frac{\big[t+\frac{\sum_{i=1}^{m-2}\gamma_{i}\zeta_{i}}{\delta}\big]}{\sigma}=\theta_{1}p(t).
\quad
\end{equation}
Therefore, $\theta_{1}p(t)\leq\mathcal{T}p(t)\leq\theta_{2}p(t),~t\in[0,1], ~that~is,~\theta_{1}p\leq\mathcal{T}p\leq\theta_{2}p, ~Hence,~\mathcal{T}p\in E_{p}.$
\par Next, for $\ell\in(0,1),~x\in E,$ and by $(\mathcal{H}_{4})$ it follows that
\begin{align*}
\mathcal{T}(\ell x)(t)&=\int_{0}^{1}G(t,q\tau;\zeta_{i})~h(\tau)~f(\tau,\ell x(\tau))~d_{q}\tau
\quad\quad\quad\quad\quad\quad\quad\quad\quad\quad\quad\\&
\geq y(\ell)\int_{0}^{1}G(t,q\tau;\zeta_{i})~h(\tau)~f(\tau, x(\tau))~d_{q}\tau\quad\quad\quad\quad\quad\quad\quad\quad\quad\quad\quad
\end{align*}
\begin{equation}
=y(\ell)\mathcal{T}(x)(t),~t\in[0,1]\quad\quad\quad\quad\quad\quad\quad\quad\quad\quad\quad\quad\quad\quad\quad
\end{equation}
So, we have $\mathcal{T}(\ell x)\geq y(\ell)\mathcal{T}x,~\forall x\in E,~\ell\in(0,1).$ Hence, all the assumptions of
Lemma 6 are satisfied. Finally, from Lemma 7, there exists a unique $x^{*}_{\lambda}\in E_{p}$ such that $\mathcal{T}x^{*}_{\lambda}=\frac{1}{\lambda}x^{*}_{\lambda},$ that is, $\lambda\mathcal{T}x^{*}_{\lambda}=x^{*}_{\lambda},$ so
\begin{equation}
x^{*}_{\lambda}(t)=\lambda\int_{0}^{1}G(t,q\tau;\zeta_{i})~h(\tau)~f(\tau, x^{*}_{\lambda}(\tau))~d_{q}\tau,~t\in[0,1].
\end{equation}
Considering Lemma 4, $x^{*}_{\lambda}$ is a unique positive solution of problem (1.1) for given $\lambda>0.$ By Lemma 7 (i), $x^{*}_{\lambda}$ is
strictly decreasing in $\lambda,$ that is, $0<\lambda_{1}<\lambda_{2}$ implies $x^{*}_{\lambda_{1}}\leq x^{*}_{\lambda_{2}},~x^{*}_{\lambda_{1}}\neq x^{*}_{\lambda_{2}}.$ And here one has $\lim_{\lambda\rightarrow+\infty} \|x^{*}_{\lambda}\|=+\infty,~\lim_{\lambda\rightarrow0^{+}} \|x^{*}_{\lambda}\|=0.$
Moreover, if there exists $r\in(0,1)$ such that $y(\ell)\geq\ell^{r}$ for $\ell\in(0,1),$ then, by using Lemma
7, it yields that $x^{*}_{\lambda}$ is continuous in $\lambda,$ that is, $\|x_{\lambda}-x_{\lambda_{0}}\|\rightarrow0~as~\lambda\rightarrow\lambda_{0}(\lambda_{0}>0).$
\par Now, let $\mathcal{T}_{\lambda}=\lambda\mathcal{T}.$ And here for $\mathcal{T}_{\lambda},$ all assumptions of Lemma 6, are satisfied.
Hence, for any initial value $x_{0}\in E_{p},$ constructing the sequence $x_{n}=\mathcal{T}_{\lambda}x_{n-1},~n=1,2,...,$ one has $x_{n}\rightarrow x^{*}_{\lambda}~as~n\rightarrow+\infty.$ Thus, the equation (3.2) holds, and $\lim_{n\rightarrow+\infty}x_{n}(t)=x^{*}_{\lambda}(t).
\quad\quad\quad\quad\quad\quad\quad\quad\quad\quad\quad\quad\quad\quad\quad\quad\quad\quad\quad\quad\quad\quad\quad\quad\quad\Box$
\\\textbf{\ Corollary 1.} Suppose that $(\mathcal{H}_{1})-(\mathcal{H}_{4})$ holds. Then the following fractional $q$-difference equation
m-point nonlinear boundary conditions
\begin{equation}
\left\{\begin{matrix} _{C}D^{\alpha}_{q}x(t)+h(t)f(t,x(t))=0,~~ t\in(0,1),~~n-1<\alpha\leq n,~~~ n>2,\\\\x(0)=\sum_{i=1}^{m-2}~\gamma_{i}~x(\zeta_{i}),\quad\quad\quad\\\\
_{C}D^{2}_{q}x(0)=~_{C}D^{3}_{q}x(0)=...=~_{C}D^{n-1}_{q}x(0)=0,\quad\quad\quad\\\\
\nu ~_{C}D_{q}x(1)-\mu\tilde{\alpha}[x]=\sum_{i=1}^{m-2} \beta_{i}~ _{C}D_{q}x(\zeta_{i}),\quad\quad\quad\end{matrix}\right.
\end{equation}
has a unique positive solution $x^{*}$ in $E_{p},$ where $p(t)=\frac{\big[t+\frac{\sum_{i=1}^{m-2}\gamma_{i}\zeta_{i}}{\delta}\big]}{\sigma},~t\in[0,1].$
Furthermore, for $x_{0}\in E_{p},$ constructing the sequence
\begin{equation}
x_{n}(t)=\int_{0}^{1}G(t,q\tau;\zeta_{i})~h(\tau)~f(\tau,x_{n-1}(\tau))~d_{q}\tau,~n=1,2,...,
\end{equation}
and here $\lim_{n\rightarrow+\infty}x_{n}(t)=x^{*}(t),~t\in[0,1].$
\\\textbf{\ Corollary 2.} Suppose that $(\mathcal{H}_{2})-(\mathcal{H}_{4})$ holds. Then the following fractional $q$-difference equation
m-point nonlinear boundary conditions
\begin{equation}
\left\{\begin{matrix} _{C}D^{\alpha}_{q}x(t)+\lambda f(t,x(t))=0,~~ t\in(0,1),~~n-1<\alpha\leq n,~~~ n>2,\\\\ _{C}D^{2}_{q}x(0)=~_{C}D^{3}_{q}x(0)=...=~_{C}D^{n-1}_{q}x(0)=0,\quad\quad\quad\quad\\\\
x(0)=\sum_{i=1}^{m-2}~\gamma_{i}~x(\zeta_{i}),\quad \nu ~_{C}D_{q}x(1)=0,
\quad\quad\quad\quad\quad\end{matrix}\right.
\end{equation}
has a unique positive solution $x^{*}$ in $E_{p},$ where $p(t)=\frac{\big[t+\frac{\sum_{i=1}^{m-2}\gamma_{i}\zeta_{i}}{\delta}\big]}{\sigma},~t\in[0,1].$
Furthermore, for $x_{0}\in E_{p},$ constructing the sequence
\begin{equation}
x_{n}(t)=\lambda\int_{0}^{1}\big[H_{1}(t,q\tau)+H_{2}(t,q\tau;\zeta_{i})\big]f(\tau,x_{n-1}(\tau))~d_{q}\tau,~n=1,2,...,
\end{equation}
and here $\lim_{n\rightarrow+\infty}x_{n}(t)=x^{*}_{\lambda}(t),~t\in[0,1].$
\\\textbf{\ Remark 2.} For the assumption $(\mathcal{H}_{4}),$ if $y(\ell)=\ell^{r},~\ell\in(0,1),$ then for any initial
value $x_{0}\in E_{p},$ there have $x_{n}\rightarrow x^{*}_{\lambda}~as~n\rightarrow\infty.$ Moreover, by a proof of Theorem 2.1 [], we have the error estimation
\[\|x_{n}-x^{*}_{\lambda}\|=o(1-\epsilon^{r^{n}}),~n\rightarrow\infty,\]
here $\epsilon\in(0,1)$ is a constant depend on $x_{0}.$
\\\textbf{\ Example 1.} Consider the following boundary value
problem for generalized Caputo-type fractional $q$-difference equation
\begin{equation}
\left\{\begin{matrix} _{C}D^{\frac{7}{2}}_{q}x(t)+\lambda h(t)f(t,x(t))=0,~~ t\in(0,1),\quad\quad\quad\quad\quad\quad\quad\quad\\\\x(0)=\frac{1}{5}x(\frac{1}{3})+\frac{1}{3}x(\frac{1}{2}),~~
 _{C}D^{2}_{q}x(0)=~_{C}D^{3}_{q}x(0)=0,\quad\quad\quad\quad\quad\\\\
 \nu ~_{C}D_{q}x(1)-\mu\tilde{\alpha}[x]=\frac{2}{3}~_{C}D_{q}x(\frac{1}{3})+\frac{1}{4}~_{C}D_{q}x(\frac{1}{2}).\quad\quad\quad\end{matrix}\right.
\end{equation}
Here, $\alpha=\frac{7}{2},q=\frac{1}{3},
\nu=5,\mu=3,\gamma_{1}=\frac{1}{5},\gamma_{2}=\frac{1}{3},\beta_{1}=\frac{2}{3},\beta_{2}=\frac{1}{4},\zeta_{1}=\frac{1}{3},\zeta_{2}=\frac{1}{2},$ and let
\begin{equation}
h(t)=\log_{e}(\frac{1}{t}),~f(t,x(t))=\bigg(\sqrt[3]{3}t^{3}+\big[x(t)\big]^{\frac{4}{3}}t+1\bigg)^{\frac{1}{3}},~t\in(0,1).
\end{equation}
Clearly,
\\$f(t,0)=\bigg(\sqrt[3]{3}t^{3}+1\bigg)^{\frac{1}{3}}>0,~h(t)>0$ and $f(t,x)$  is increasing in $x$ for $x\in\mathbb{R}^{+},~t\in[0,1].$ Therefore, assumptions
$(\mathcal{H}_{2}),(\mathcal{H}_{3})~$hold. Set $y(\ell)=\ell^{\frac{4}{9}},$ we have $y(\ell)\in(\ell,1),~\ell\in(0,1).$ Thus, for $x\in[0,\infty),$ we obtain
\begin{equation}
f(t,\ell x)=\bigg(\sqrt[3]{3}t^{3}+\big[\ell x(t)\big]^{\frac{4}{3}}t+1\bigg)^{\frac{1}{3}}\geq\ell^{\frac{4}{9}}\bigg(\sqrt[3]{3}t^{3}+\big[\ell x(t)\big]^{\frac{4}{3}}t+1\bigg)^{\frac{1}{3}}=y(\ell)f(t,x).
\end{equation}
For the Riemann-Stieltjes integral we can discuss two cases:
\\$(1)$ Let $\tilde{\alpha}[x]=0,$ then we get $\mathbf{B}=0, \nu,\mu>0,$ and $\rho=\nu-\sum_{i=1}^{2} \beta_{i}=\frac{49}{12}>0.$
\\$(2)$ Let $\tilde{\alpha}[x]=\int_{0}^{1}[2t-1]x(t)dt,$ here the sign of the function $(2t-1)$ is change for $t\in[0,1].$ Therefore, we have \\$\mathbf{B}=\int_{0}^{1}\big[t+\frac{\sum_{i=1}^{m-2}~\gamma_{i}\zeta_{i}}{\delta}\big]\big(2t-1\big)dt=\int_{0}^{1}\big(2t^{2}-\frac{1}{2}\big)dt=\frac{1}{6}>0,$ and $\rho\simeq3.5833>0.$
All the assumptions of Theorem 1 are satisfied, then the problem (3.12) has a unique positive solution $x^{*}_{\lambda}\in  E_{p},$ and here $p(t)=\frac{1}{3}(2t+1).$ Moreover, for $x_{0}\in  E_{p},$ we can construct
the sequence
\[x_{n}(t)=\lambda\int_{0}^{1}G(t,q\tau;\zeta_{i})~\log_{e}(\frac{1}{\tau})~
\big(\sqrt[3]{3}\tau^{3}+\big[x_{n-1}(\tau)\big]^{\frac{4}{3}}\tau+1\big)^{\frac{1}{3}}~d_{q}\tau,~n=1,2,...,\]
so, $x_{n}(t)\rightarrow x^{*}_{\lambda}(t),~as~n\rightarrow\infty.$

  \end{document}